\documentclass[10pt]{article}
\usepackage[T1]{fontenc}
\usepackage{lmodern}
\usepackage{amsmath,amssymb,amsthm,mathtools}
\usepackage{microtype}
\usepackage[a4paper,margin=1in]{geometry}
\usepackage{enumitem}
\usepackage{needspace}
\usepackage{float}
\usepackage{tikz}
\usetikzlibrary{arrows.meta,fit,backgrounds}
\usepackage{hyperref}
\hypersetup{hidelinks,pdfauthor={Xuanying Han and Sophie C.C. Sun}}
\setlist{nosep}

\numberwithin{equation}{section}
\newtheorem{theorem}{Theorem}[section]
\newtheorem{lemma}[theorem]{Lemma}
\newtheorem{proposition}[theorem]{Proposition}
\newtheorem{corollary}[theorem]{Corollary}
\theoremstyle{definition}
\newtheorem{definition}[theorem]{Definition}
\newtheorem{example}[theorem]{Example}
\theoremstyle{remark}
\newtheorem{remark}[theorem]{Remark}

\newcommand{\Fn}{\mathcal{F}_n}
\newcommand{\Anj}{\mathcal{A}_{n,j}}
\newcommand{\divw}{\operatorname{div}_w}

\title{\bfseries UNDERLYING TREES OF FLOOR DIAGRAMS IN~GENUS~ZERO}
\author{Xuanying Han and Sophie C.C. Sun}
\date{}

\begin{document}
\maketitle
\vspace{-2.3em}

\begin{abstract}
We count the distinct trees on the vertex set $\{1,\ldots,n\}$ obtained by forgetting the edge weights of labeled floor diagrams of genus zero. Answering a question of Fomin and Mikhalkin, we prove that these trees are equinumerous with the alternating trees on the same vertex set. Their common number is $\frac{1}{n\,2^{n-1}}\sum_{k=1}^{n}\binom{n}{k}k^{n-1}$.
By decomposing the underlying trees at a root, we derive a differential equation that determines their generating function by number of vertices. Extracting its coefficients gives the counting formula above. We also obtain an explicit expression for the generating function that refines this count by root rank and capacity.
\end{abstract}

\medskip

\section{Introduction}
Mikhalkin's correspondence theorem~\cite{Mik05} relates algebraic counts of plane curves to tropical ones. Floor diagrams express these tropical counts in terms of weighted directed graphs~\cite{BM,BM09,FM}. They have been used to compute plane node polynomials~\cite{Block11} and to prove eventual polynomiality of Severi degrees for families of toric surfaces~\cite{AB13}. In genus zero the underlying graph is a tree. Fomin and Mikhalkin~\cite[Theorem~6.1]{FM} showed that there are $n^{n-2}$ labeled floor diagrams of degree $n$ and genus zero, with a single diagram for $n=1$. Different edge weightings can give the same labeled tree. We ask how many distinct trees remain when the weights are forgotten.

For a tree $T$ with vertex set $V(T)=\{1,\ldots,n\}$, orient each edge from its smaller endpoint to its larger endpoint.
Assign a positive integer weight $w(e)$ to each edge $e$ of $T$. The divergence at a vertex $v$ is defined by
\begin{equation}\label{eq:div}
 \divw(v)=\sum_{v\to u} w(vu)-\sum_{u\to v} w(uv).
\end{equation}
A weighting $w$ is \emph{feasible} if every vertex $v$ of $T$ satisfies
\begin{equation}\label{eq:feasible}
 \divw(v)\le 1.
\end{equation}
Let $\Fn$ consist of all trees on $\{1,\ldots,n\}$ that admit a feasible weighting, and set $f_n=|\Fn|$. By~\cite[Definition~1.1]{FM}, these are the underlying trees of labeled floor diagrams of degree $n$ and genus zero.

An alternating tree is a labeled tree in which the neighbors of each vertex are either all smaller or all larger than that vertex.

\begin{figure}[H]
\centering
\begin{tikzpicture}[
  vertex/.style={circle,draw,minimum size=6mm,inner sep=0pt,font=\small},
  edge/.style={-{Stealth[length=1.8mm,width=1.4mm]},line width=.55pt}
]
  \node[vertex] (a1) at (0,0) {$1$};
  \node[vertex] (a2) at (1.3,0) {$2$};
  \node[vertex] (a3) at (2.6,0) {$3$};
  \draw[edge] (a1) -- (a2);
  \draw[edge] (a2) -- (a3);
  \node at (1.3,-1.2) {$T_1$};

  \node[vertex] (b1) at (5,.7) {$1$};
  \node[vertex] (b2) at (5,-.7) {$2$};
  \node[vertex] (b3) at (6.8,0) {$3$};
  \draw[edge] (b1) -- (b3);
  \draw[edge] (b2) -- (b3);
  \node at (5.9,-1.2) {$T_2$};

  \node[vertex] (c1) at (9.2,0) {$1$};
  \node[vertex] (c2) at (11,.7) {$2$};
  \node[vertex] (c3) at (11,-.7) {$3$};
  \draw[edge] (c1) -- (c2);
  \draw[edge] (c1) -- (c3);
  \node at (10.1,-1.2) {$T_3$};
\end{tikzpicture}

\caption{The three trees on $\{1,2,3\}$.}
\label{fig:capacity-n3}
\end{figure}
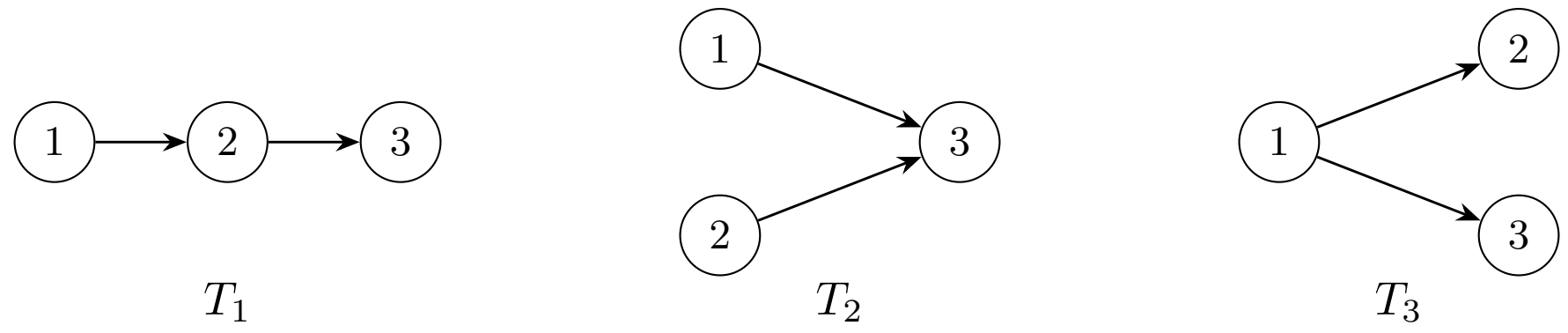

Figure~\ref{fig:capacity-n3} shows the three trees on $\{1,2,3\}$. On the path $T_1$, the edge weights $(1,1)$ and $(1,2)$, in path order, are both feasible. Unit weights also make $T_2$ feasible. In $T_3$ both edges leave vertex $1$, so its divergence is at least $2$ for every choice of weights.

Thus $\mathcal F_3=\{T_1,T_2\}$, while the alternating trees on $\{1,2,3\}$ are $T_2$ and $T_3$. Fomin and Mikhalkin~\cite[Problem~7.4]{FM} asked whether these two families have the same size for every $n$.

Postnikov~\cite[Theorem~1]{Postnikov} obtained an explicit formula for the number of alternating trees, and Chauve, Dulucq, and Rechnitzer~\cite{CDR01} gave a bijective proof of that formula. The following theorem shows that the same formula counts the distinct underlying trees of labeled floor diagrams of genus zero.

\begin{theorem}\label{thm:main}
For every $n\ge1$,
\begin{equation}\label{eq:fn}
 f_n=\frac{1}{n2^{n-1}}\sum_{k=1}^n \binom{n}{k}k^{n-1}.
\end{equation}
\begin{samepage}
Equivalently, the formal power series
\begin{equation}\label{eq:Adef}
 A(z)=\sum_{n\ge1} f_n\frac{z^{n-1}}{(n-1)!}
\end{equation}
satisfies
\begin{equation}\label{eq:Aeq}
 A(z)=\exp\!\left(\frac{z(1+A(z))}{2}\right).
\end{equation}
In particular, $\Fn$ and the alternating trees with vertex labels $1,\ldots,n$ are equinumerous.
\end{samepage}
\end{theorem}

To count the trees in $\Fn$ recursively, we need to know which weights
a branch allows on the edge joining it to the root. In
Section~\ref{sec:feasibility} we use the charge of an admissible rooted
tree to describe these weights and express the charge of a tree in terms
of the charges of its branches. In Section~\ref{sec:root-gf} we track
the rank of the root and apply the exponential formula for labeled sets
to obtain an integral identity, which we differentiate to derive a
differential equation.
In Section~\ref{sec:main-proof} we use this differential equation to recover Postnikov's functional equation~\eqref{eq:Aeq} for alternating trees~\cite[Theorem~2]{Postnikov} and extract the coefficients of $A(z)$ to complete the proof of Theorem~\ref{thm:main}. Theorem~\ref{thm:refinement} then gives explicit formulas for the three series that determine the generating function by root rank and capacity.

For $n\ge2$, $f_n$ also counts the regions of the Linial arrangement on $n-1$ coordinates, defined by the hyperplanes $x_i-x_j=1$ for $1\le i<j\le n-1$. This follows from Theorem~\ref{thm:main} and the region count in~\cite[Theorem~8.2]{PS}. For related work on this arrangement, see~\cite{Ath96,Bernardi18,Stanley07}.

Throughout the paper, all generating functions and differential equations are understood formally.

\section{A feasibility recursion on rooted trees}
\label{sec:feasibility}
We seek a recursive criterion for feasibility in terms of the branches
at a chosen root.
Since deleting the root of a feasible tree can leave an infeasible branch,
we introduce a weaker condition on rooted trees.

\begin{definition}\label{def:charge}
A rooted labeled tree $(T,r)$ is \emph{admissible} if there is a positive
integral edge weighting $w$ such that $\divw(v)\le1$ for every $v\ne r$.
For an admissible $(T,r)$, its \emph{charge} is
\begin{equation}\label{eq:qdef}
 q(T,r)=1-\min_w \divw(r),
\end{equation}
where the minimum is taken over all such weightings.

If $T\in\Fn$ and $r$ is a vertex of $T$, the \emph{root capacity} of $T$ at $r$ is
\[
 \kappa_T(r)=
 1-\min_{w\text{ feasible on }T}\divw(r).
\]
\end{definition}

For a tree on $n$ vertices, the edge contributions cancel in the sum of
all divergences. For any weighting satisfying the bounds away from $r$,
this gives
\[
 \divw(r)=-\sum_{v\ne r}\divw(v)\ge -(n-1).
\]
For an admissible rooted tree, the possible root divergences form a nonempty set of integers. The lower bound ensures that the minimum in~\eqref{eq:qdef} is attained and gives $q(T,r)\le n$. For feasible weightings, the bound at $r$ gives $-(n-1)\le\divw(r)\le1$. In particular, $0\le\kappa_T(r)\le n$ for $T\in\Fn$.

\begin{remark}\label{lem:feasibility}
For every choice of root $r$, a labeled tree $T$ belongs to $\Fn$
if and only if $(T,r)$ is admissible and $q(T,r)\ge0$. In this case,
\[
  q(T,r)=\kappa_T(r).
\]
\end{remark}

Deleting $r$ leaves branches $T_i$, each rooted at the vertex $u_i$
that was adjacent to $r$. Each branch inherits the relative order of
its vertex labels and hence its edge orientations from $T$, so the
definitions of admissibility and charge apply with $u_i$ as root.

\begin{lemma}\label{lem:branch}
The rooted tree $(T,r)$ is admissible if and only if each branch $(T_i,u_i)$ is admissible and
\begin{equation}\label{eq:branch-cond}
 q(T_i,u_i)\ge1\qquad\text{whenever }u_i<r.
\end{equation}
When these conditions hold, writing $q_i=q(T_i,u_i)$ gives
\begin{equation}\label{eq:branch-rec}
 q(T,r)=1+\sum_{u_i<r}q_i-\sum_{u_i>r}\max(1,-q_i).
\end{equation}
For the tree with one vertex, the charge is $1$.
\end{lemma}

\begin{proof}
Attaching $T_i$ to $r$ changes the divergence at $u_i$ and nowhere
else in $T_i$. Thus, if $(T,r)$ is admissible, each rooted branch is
admissible. Now suppose that all branches are admissible, and write
$q_i=q(T_i,u_i)$. The least attainable divergence at $u_i$ within $T_i$
is $1-q_i$.

If $u_i<r$, the edge $u_i\to r$ adds its weight $s$ to the divergence at $u_i$. The bound there can be met exactly when $(1-q_i)+s\le1$, so the allowable positive weights are the integers $1\le s\le q_i$. Such a weight exists precisely when $q_i\ge1$.

If $u_i>r$, the edge $r\to u_i$ subtracts $s$ from the divergence at $u_i$. In this case the bound can be met exactly when $(1-q_i)-s\le1$, and the allowable weights are the integers $s\ge\max(1,-q_i)$.

To minimize the divergence at $r$, choose the largest allowed weight $q_i$ on each edge entering $r$ and the smallest allowed weight $\max(1,-q_i)$ on each edge leaving $r$. The branches have disjoint vertex sets, so all these extrema can be achieved simultaneously. Substituting the resulting minimum in~\eqref{eq:qdef} gives~\eqref{eq:branch-rec}.
\end{proof}

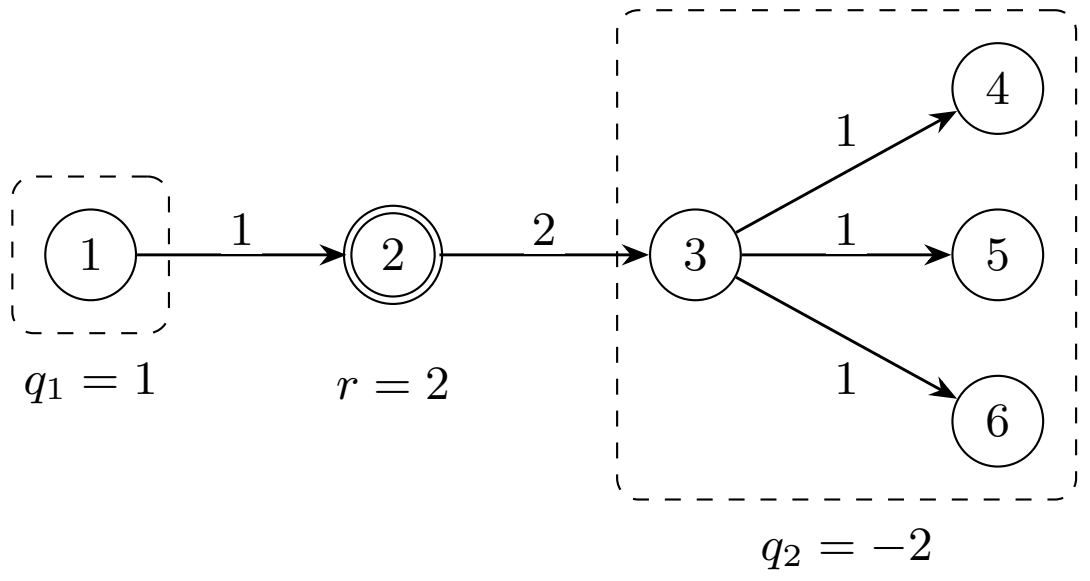
\begin{figure}[H]
\centering
\begin{tikzpicture}[
  vertex/.style={circle,draw,minimum size=6mm,inner sep=0pt,font=\small},
  edge/.style={-{Stealth[length=1.8mm,width=1.4mm]},line width=.55pt},
  weight/.style={font=\small,fill=white,inner sep=1.5pt},
  branch/.style={draw,dashed,rounded corners=4pt,inner sep=6pt}
]
  \node[vertex] (v1) at (0,0) {$1$};
  \node[vertex,double,double distance=1pt] (v2) at (2,0) {$2$};
  \node[vertex] (v3) at (4,0) {$3$};
  \node[vertex] (v4) at (6,1.1) {$4$};
  \node[vertex] (v5) at (6,0) {$5$};
  \node[vertex] (v6) at (6,-1.1) {$6$};
  \draw[edge] (v1) -- node[weight,above] {$1$} (v2);
  \draw[edge] (v2) -- node[weight,above] {$2$} (v3);
  \draw[edge] (v3) -- node[weight,midway,above=3pt] {$1$} (v4);
  \draw[edge] (v3) -- node[weight,above] {$1$} (v5);
  \draw[edge] (v3) -- node[weight,midway,below=3pt] {$1$} (v6);
  \begin{scope}[on background layer]
    \node[branch,fit=(v1)] {};
    \node[branch,fit=(v3)(v4)(v5)(v6)] {};
  \end{scope}
  \node[font=\small] at (0,-.85) {$q_1=1$};
  \node[font=\small] at (2,-.85) {$r=2$};
  \node[font=\small] at (5,-1.95) {$q_2=-2$};
\end{tikzpicture}
\caption{The two branches at root $2$ have charges $1$ and $-2$.}
\label{fig:branch-example}
\end{figure}

\begin{example}\label{ex:unit}
Consider the tree $T$ on six vertices in Figure~\ref{fig:branch-example},
rooted at vertex $2$. Deleting the root leaves the two branches shown in
dashed boxes. One is the singleton $1$, of charge $q_1=1$, and the other
is a star rooted at $3$.

All three edges of the star point away from $3$, so the minimum internal divergence at $3$ is $3$, attained when all three weights are $1$. The star therefore has charge $q_2=1-3=-2$. To satisfy the divergence bound at vertex $3$ when the star is reattached, the incoming edge $2\to3$ must carry weight at least $2$. By~\eqref{eq:branch-rec},
\[
 q(T,2)=1+1-2=0.
\]
Taking weight $2$ on $2\to3$ and weight $1$ on every other edge attains this charge. The divergences in label order are $1,1,1,-1,-1,-1$, so $T\in\mathcal F_6$.
\end{example}

\section{A generating function for root rank and charge}
\label{sec:root-gf}
We use the charge recursion to derive an integral identity for admissible
trees counted by root rank and charge, and from it a differential
equation for the generating function of feasible trees.

\subsection{The exponential formula}
In~\eqref{eq:branch-rec}, the contribution of a branch depends on whether
its root label is smaller or larger than that of the deleted root. The
\emph{rank} of the root is the position of its label in the increasing
order of all $n$ vertex labels. Thus, for a root of rank $j$, exactly
$j-1$ labels are smaller than its label and $n-j$ are larger.
On the standard vertex set $\{1,\ldots,n\}$, vertex $j$ has rank $j$.
Let $\Anj$ be the set of admissible trees on this vertex set rooted at $j$.

Assign real labels in $[0,1]$ to the vertices and fix the root label at $x$.
For root rank $j$, we use the volume of the region realizing one prescribed
strict ordering of the labels in which $x$ occupies position $j$.
With $t$ marking the charge, define
\begin{equation}\label{eq:Fdef}
 F(x,t;z)=\sum_{n\ge1}z^n\sum_{j=1}^n
 \frac{x^{j-1}(1-x)^{n-j}}{(j-1)!(n-j)!}
 \sum_{T\in\Anj} t^{q(T,j)}.
\end{equation}
For each $n$, $[z^n]F$ is a polynomial in $x$ and a Laurent polynomial in $t$,
so $F\in z\mathbb{Q}[x,t,t^{-1}][[z]]$. Split $F$ according to the sign of
the charge, writing
\begin{equation}\label{eq:FGparts}
 F=F_++F_0+F_-,\qquad G=F_++F_0,
\end{equation}
where the subscripts refer to positive, zero, and negative powers of $t$. By Remark~\ref{lem:feasibility}, $G$ counts feasible trees, with $t$ recording the root capacity. Write $a(x)=F_0(x;z)=G(x,0;z)$ for the series counting feasible trees whose root has capacity zero. We suppress its dependence on $z$ in the notation $a(x)$.

With $x$ fixed, let the remaining labels
$y_1,\ldots,y_{n-1}$ vary in $[0,1]$. For a prescribed strict ordering
in which $x$ has rank $j$, the labels below and above $x$ range
independently over two ordered simplices.
Their volumes are $\frac{x^{j-1}}{(j-1)!}$ and $\frac{(1-x)^{n-j}}{(n-j)!}$,
respectively. Their product gives the rank factor in~\eqref{eq:Fdef}.
At $x=0$ only rank $1$ occurs, and at $x=1$ only rank $n$ occurs.

For $0<x<1$ and $n=j=3$, both $y_1$ and $y_2$ lie below $x$. Their
two possible orders give the two triangles in
Figure~\ref{fig:volume-example}(a), each of area $\frac{x^2}{2}$. For $n=3$
and $j=1$, both labels lie above $x$. The diagonal of $[x,1]^2$
separates the two possible orders, $x<y_1<y_2$ and $x<y_2<y_1$,
giving the two triangles in Figure~\ref{fig:volume-example}(b),
each of area $\frac{(1-x)^2}{2}$.

\begin{figure}[H]
\centering
\begin{tikzpicture}[>=Latex,guide/.style={gray!55,line width=.4pt,dash pattern=on 2.5pt off 2pt}]
  \def\unitend{4.6}
  \newcommand{\orderaxes}{%
    \draw[guide] (0,\unitend) -- (\unitend,\unitend) -- (\unitend,0);
    \draw[->] (-.08,0) -- (5.02,0) node[right] {$y_1$};
    \draw[->] (0,-.08) -- (0,5.02) node[above] {$y_2$};
    \node[below left] at (0,0) {$0$};
    \foreach \pos/\lab in {\rootpos/x,\unitend/1}{
      \draw (\pos,0) -- ++(0,-.065) node[below] {$\lab$};
      \draw (0,\pos) -- ++(-.065,0) node[left] {$\lab$};
    }
  }
  \begin{scope}[scale=.8]
    \def\rootpos{3.45}
    \fill[gray!22] (0,0) -- (0,\rootpos) -- (\rootpos,\rootpos) -- cycle;
    \orderaxes
    \draw[thick] (0,0) rectangle (\rootpos,\rootpos);
    \draw[thick] (0,0) -- (\rootpos,\rootpos);
    \node at (1.15,2.30) {$y_1<y_2$};
    \node at (2.30,-.78) {(a)\enspace $j=3$};
  \end{scope}
  \begin{scope}[xshift=7cm,scale=.8]
    \def\rootpos{.92}
    \fill[gray!22] (\rootpos,\rootpos) -- (\rootpos,\unitend)
                  -- (\unitend,\unitend) -- cycle;
    \orderaxes
    \draw[guide] (0,\rootpos) -- (\rootpos,\rootpos)
                                    -- (\rootpos,0);
    \draw[thick] (\rootpos,\rootpos) rectangle (\unitend,\unitend);
    \draw[thick] (\rootpos,\rootpos) -- (\unitend,\unitend);
    \node at (2.147,3.373) {$y_1<y_2$};
    \node at (2.30,-.78) {(b)\enspace $j=1$};
  \end{scope}
\end{tikzpicture}
\caption{Order regions when the root has the largest label (a) or the smallest label (b).}
\label{fig:volume-example}
\end{figure}

Let $m=n-1$. Given a tree $T$ on $\{1,\ldots,m+1\}$ rooted at $j$,
give its root index $0$ and assign the indices $1,\ldots,m$ bijectively
to the remaining vertices. This gives a tree $H$ on $\{0,\ldots,m\}$
rooted at $0$. Record the indices in increasing order of the original
vertex labels. To reverse the construction, give the $k$th vertex in
this order the label $k$ and orient each edge from its smaller label
to its larger label. This recovers $T$ together with its chosen indexing.
Thus trees $T$ rooted at $j$ with a choice of indexing are in bijection
with pairs consisting of $H$ and an order on its vertices in which $0$
has rank $j$. For a general treatment of structures on totally ordered
sets, see Bergeron, Labelle, and Leroux~\cite[Chapter~5]{BLL}.

We now combine this correspondence with the volumes computed above to
rewrite~\eqref{eq:Fdef} as an integral. Give vertex $0$ of $H$ the real
label $x$ and vertex $i$ the real label $y_i$ for $1\le i\le m$.
For distinct labels, the inverse construction recovers $T$ with root
label $j$ equal to the rank of $x$.
For $\boldsymbol{y}=(y_1,\ldots,y_m)$, set
\[
 J_H(x,\boldsymbol{y};t)=
 \begin{cases}
 t^{q(T,j)}, & \text{if the labels are distinct and $(T,j)$ is admissible},\\
 0,          & \text{otherwise}.
 \end{cases}
\]
\Needspace{6\baselineskip}
For a fixed $x$ and an attainable root rank $j$, the $m!$ indexings of
a given tree $T\in\mathcal A_{m+1,j}$ give $m!$ distinct orders of the
indices. The region in $[0,1]^m$ for each order has volume equal to the factor
for rank $j$ in~\eqref{eq:Fdef}, and $J_H$ equals $t^{q(T,j)}$
throughout the region. Summing over $H$ and integrating over
$\boldsymbol{y}$ counts one contribution from each indexing of $T$.
Dividing by $m!$ therefore gives the contribution of $T$ to
$[z^{m+1}]F$. Summing over $m$ gives the coefficientwise integral
representation
\begin{equation}\label{eq:Fint}
 F(x,t;z)=z\sum_{m\ge0}\frac{z^m}{m!}
   \sum_H\int_{[0,1]^m}J_H(x,\boldsymbol{y};t)\,d\boldsymbol{y}.
\end{equation}

Deleting the root leaves an unordered set of branches, so we can apply the exponential formula for labeled sets~\cite[Chapter~II]{FS}.

\Needspace{12\baselineskip}
\begin{lemma}\label{lem:exp}
The function $F$ satisfies
\begin{equation}\label{eq:expformula}
 F(x,t;z)=zt\exp\!\left(
 \int_0^x F_+(y,t;z)\,dy
 +\int_x^1\left(F_-(y,t;z)+\frac{G(y,1;z)}{t}\right)\,dy
 \right).
\end{equation}
\end{lemma}

\begin{proof}
Fix the real label $x$ of the root $r$. Let $B(x,t;z)$ be the generating
function for one branch attached to $r$. By Lemma~\ref{lem:branch}, a
branch with root label $y<x$ must have positive charge, so these branches
contribute $F_+(y,t;z)$. For $y>x$, a branch of negative charge $q$
retains its weight $t^q$, so these branches contribute $F_-(y,t;z)$.
In this range, every branch of nonnegative charge contributes the same
factor $\frac{1}{t}$, giving $\frac{G(y,1;z)}{t}$ in total. Integrating over
the two ranges of $y$ gives
\[
 B(x,t;z)=\int_0^x F_+(y,t;z)\,dy+\int_x^1\left(F_-(y,t;z)+\frac{G(y,1;z)}{t}\right)\,dy.
\]
Write $B_s(x,t)=[z^s]B(x,t;z)$ for the contribution from branches of
size $s$. We abbreviate $B_s(x,t)$ to $B_s$ in the products below.

Temporarily order the $p\ge1$ branches and fix their positive sizes
$s_1,\ldots,s_p$, with $s_1+\cdots+s_p=m$. There are
$\binom{m}{s_1,\ldots,s_p}$ ways to distribute the $m$ nonroot vertices
among the branches, and branch $i$ has $s_i$ possible roots. Including
the factor $\frac{1}{m!}$ from~\eqref{eq:Fint} gives
\[
 \frac{1}{m!}\binom{m}{s_1,\ldots,s_p}\prod_{i=1}^p s_i
 =\prod_{i=1}^p\frac{s_i}{s_i!}
 =\prod_{i=1}^p\frac{1}{(s_i-1)!}.
\]
For branch $i$, the factor $\frac{1}{(s_i-1)!}$ is precisely the normalization
in~\eqref{eq:Fint} with its root fixed. In $B_{s_i}$, the label of that
root is integrated over the appropriate interval. For fixed $x$, the
integral factors over the branches. Let $T$ be the tree obtained by
attaching these branches to $r$. By~\eqref{eq:branch-rec}, the branch
contributions to the charge sum to $q(T,r)-1$, so their powers of $t$
multiply to $t^{q(T,r)-1}$. After integration, the sum over trees on
each rooted vertex set depends only on its size. Thus the ordered branches
contribute $\prod_{i=1}^p B_{s_i}$.

Summing over all ordered lists of positive sizes, with $z$ marking the
total size, gives $B(x,t;z)^p$. Each unordered family of $p$ branches
has exactly $p!$ orderings, since the branches have distinct vertex sets.
Dividing by $p!$, restoring the root factor $zt$, and summing over $p$
gives
\[
 F(x,t;z)=zt\sum_{p\ge0}\frac{B(x,t;z)^p}{p!}
         =zt\exp\!\bigl(B(x,t;z)\bigr).
\]
The term $p=0$ accounts for the tree consisting only of the root.
Substituting the integral expression for $B$ gives~\eqref{eq:expformula}.
\end{proof}

\subsection{Feasible roots}
Setting $t=1$ in $G$ sums over all capacities at each root rank.
The resulting series is independent of the root label $x$, a fact
we use to derive a differential equation for $G$.
\begin{lemma}\label{lem:C}
The series $G(x,1;z)$ is independent of $x$. Its common value $C(z)$ is given by
\begin{equation}\label{eq:Cdef}
 G(x,1;z)=C(z)=\sum_{n\ge1}f_n\frac{z^n}{(n-1)!}=zA(z).
\end{equation}
\end{lemma}

\begin{proof}
With the root fixed at vertex $j$, Remark~\ref{lem:feasibility} identifies the
trees of nonnegative charge with the $f_n$ trees in $\Fn$.
Summing the rank factors in~\eqref{eq:Fdef} over $j$ gives
\[
 \sum_{j=1}^n \frac{x^{j-1}(1-x)^{n-j}}{(j-1)!(n-j)!}
 =\frac{(x+(1-x))^{n-1}}{(n-1)!}=\frac{1}{(n-1)!}.
\]
Multiplying by $f_n z^n$ and summing over $n$ gives~\eqref{eq:Cdef}.
\end{proof}

\begin{proposition}\label{prop:Gode}
With $a(x)=G(x,0;z)$, we have
\begin{equation}\label{eq:Gode}
 \frac{\partial G}{\partial x}=(G-a)\left(G-\frac{C}{t}\right),
\end{equation}
with the boundary conditions
\begin{equation}\label{eq:Gbc}
 G(0,t;z)=C-z+zt,\qquad a(1)=0.
\end{equation}
\end{proposition}

\begin{proof}
Differentiating~\eqref{eq:expformula} and using Lemma~\ref{lem:C} gives
\[
 \frac{\partial F}{\partial x}=F\left(F_+-F_--\frac{C}{t}\right).
\]
Substituting $F=F_++a+F_-$ and expanding gives
\[
 \frac{\partial F}{\partial x}
 =F_+^2-F_-^2+a(F_+-F_-)-\frac{C}{t}(F_++a+F_-).
\]
Since $F_+$ contains only positive powers of $t$, $\frac{F_+}{t}$ has no negative
powers. Taking the nonnegative part of the expansion gives
\[
 \frac{\partial G}{\partial x}
 =F_+^2+aF_+-\frac{C}{t}F_+
 =(G-a)\left(G-\frac{C}{t}\right).
\]

At $x=0$ the root is smallest, so each branch lowers the root charge by at
least one. Only the singleton has positive charge, contributing $zt$.
Since $G(0,1;z)=C$, we have $G(0,t;z)=C-z+zt$. At $x=1$ the root is
largest, and every branch that can be attached has positive charge
by~\eqref{eq:branch-cond}. The root charge is then at least one,
so $a(1)=0$.
\end{proof}

\section{Proof of Theorem~\ref{thm:main}}
\label{sec:main-proof}
We first regard $C=zA$ as fixed. Since both factors on the right
of~\eqref{eq:Gode} are $O(z)$, the initial value at $x=0$ determines
$G$ uniquely, one coefficient of $z$ at a time. The rational form of
$G$ in Lemma~\ref{lem:rational} then allows us to use the remaining
boundary condition $a(1)=0$ to derive a functional equation for $A$.

\begin{lemma}\label{lem:rational}
With $a(x)=G(x,0;z)$ as above, there are formal series
$b,d\in z\mathbb{Q}[x][[z]]$ such that
\begin{equation}\label{eq:rationalG}
 G(x,t;z)=a(x)+\frac{b(x)t}{1-d(x)t}.
\end{equation}
We suppress the dependence of $a,b,d$ on $z$ and use primes for
differentiation with respect to $x$ at fixed $z$. These series satisfy
\begin{equation}\label{eq:abd}
 a'=-Cb,\qquad b'=b(a-Cd),\qquad d'=b,
\end{equation}
with
\begin{equation}\label{eq:abd-init}
 a(0)=C-z,\qquad b(0)=z,\qquad d(0)=0.
\end{equation}
\end{lemma}

\begin{proof}
Seek a solution of~\eqref{eq:Gode} of the form
\[
 R=a+\frac{bt}{1-dt},
 \qquad a,b,d\in z\mathbb{Q}[x][[z]].
\]
For this candidate series, the first factor on the right
of~\eqref{eq:Gode} is $R-R(x,0;z)=R-a$. Differentiating $R$ with
respect to $x$ gives
\[
 R'=a'+\frac{b't}{1-dt}+\frac{bd't^2}{(1-dt)^2}.
\]
Substituting $R$ into the right side of~\eqref{eq:Gode} gives
\[
 (R-a)\left(R-\frac{C}{t}\right)
 =-Cb+\frac{b(a-Cd)t}{1-dt}+\frac{b^2t^2}{(1-dt)^2}.
\]
Comparing coefficients in ascending powers of $t$ gives
\[
 a'=-Cb,\qquad b'=b(a-Cd),\qquad bd'=b^2.
\]
At $x=0$, condition~\eqref{eq:Gbc} requires $R(0,t;z)=C-z+zt$.
Comparing the coefficients of $1$, $t$, and $t^2$ gives
\[
 a(0)=C-z,\qquad b(0)=z,\qquad b(0)d(0)=0.
\]
Since $b(0)=z$, the last identity gives $d(0)=0$. It also follows
that $b$ is nonzero, so $bd'=b^2$ gives $d'=b$.

We now construct $a,b,d$ from~\eqref{eq:abd} and~\eqref{eq:abd-init},
one degree in $z$ at a time. At each degree, the first two equations
determine the derivatives of the coefficients of $a$ and $b$ from lower
degrees, and the initial values fix the constants of integration.
Integrating $d'=b$ then gives the coefficient of $d$ at the same degree.
Since $d=O(z)$, the series $R$ is well defined. It satisfies
the differential equation~\eqref{eq:Gode} with the prescribed initial
value $R(0,t;z)=C-z+zt$. Uniqueness gives $R=G$, and setting $t=0$
identifies $a(x)$ with $G(x,0;z)$.
\end{proof}

\Needspace{7\baselineskip}
\begin{corollary}\label{cor:invariants}
With $a,b,d$ as in Lemma~\ref{lem:rational},
\begin{equation}\label{eq:invariants}
 a+Cd=C-z,\qquad b=(1-d)(z+Cd).
\end{equation}
\end{corollary}

\begin{proof}
Since $C$ is independent of $x$, equations~\eqref{eq:abd} give $a'+Cd'=0$.
The initial conditions~\eqref{eq:abd-init} therefore yield $a+Cd=C-z$,
or equivalently $a=C-z-Cd$.

Substituting the resulting expression for $a$ into $b'=b(a-Cd)$
in~\eqref{eq:abd} gives
\[
 b'=b(C-z-2Cd).
\]
Since $d'=b$, we replace $b$ on the right by $d'$ and expand to obtain
\[
 b'=(C-z)d'-2Cd\,d'.
\]
The two terms on the right are the derivatives of $(C-z)d$ and $-Cd^2$,
respectively. Integrating from $0$ to $x$ and using $d(0)=0$ gives
\[
 b(x)-b(0)=(C-z)d(x)-Cd(x)^2.
\]
Finally, substituting the initial value $b(0)=z$ from~\eqref{eq:abd-init}
and factoring gives the second identity
\[
 b=z+(C-z)d-Cd^2=(1-d)(z+Cd).
\]
\end{proof}

\Needspace{8\baselineskip}
\begin{proof}[Proof of Theorem~\ref{thm:main}]
Equations~\eqref{eq:abd} and~\eqref{eq:invariants} give
$d'=(1-d)(z+Cd)$. Since the denominator has constant term one, we may define
\[
 \rho(x)=\frac{z+Cd(x)}{1-d(x)}.
\]
The quotient rule and the differential equation for $d$ give
\[
 \rho'(x)=\frac{(C+z)d'(x)}{(1-d(x))^2}=(C+z)\rho(x).
\]
Since $d(0)=0$, we have $\rho(0)=z$ and hence
\begin{equation}\label{eq:rho}
 \rho(x)=z\exp((C+z)x).
\end{equation}
At $x=1$, the identity $a+Cd=C-z$ and the boundary condition $a(1)=0$
give $z+Cd(1)=C$ and $C(1-d(1))=z$. Since $C=zA$, we have
$A(1-d(1))=1$, and thus $\rho(1)=CA=zA^2$. Equation~\eqref{eq:rho}
gives a second expression for the same value, $\rho(1)=z\exp(C+z)$.
Comparing the two yields
\[
 A^2=\exp\bigl(z(1+A)\bigr).
\]
As $A(0)=1$, taking formal logarithms gives $2\log A=z(1+A)$, or
equivalently~\eqref{eq:Aeq}. For each $m\ge1$, the coefficient of $z^m$
on the right of~\eqref{eq:Aeq} depends only on the coefficients of $A$
of lower degree, so this equation determines $A$ uniquely.

To extract the coefficients, use Postnikov's
substitution~\cite[Section~2]{Postnikov},
\[
 W=z+C=z(1+A), \qquad W=z\left(1+\exp\!\left(\frac{W}{2}\right)\right).
\]
Lagrange inversion~\cite[Theorem~2.1.1]{Gessel16} gives, for $n\ge2$,
\begin{align*}
 \frac{f_n}{(n-1)!}
 &= [z^n]C=[z^n]W\\
 &=\frac1n[u^{n-1}]\left(1+\exp\!\left(\frac{u}{2}\right)\right)^n\\
 &=\frac{1}{n2^{n-1}(n-1)!}\sum_{k=1}^n\binom{n}{k}k^{n-1}.
\end{align*}
The singleton case gives $f_1=1$, also in agreement with~\eqref{eq:fn}.
Postnikov~\cite[Theorem~1]{Postnikov} obtained the same formula for the number of alternating trees with vertex labels $1,\ldots,n$.
\end{proof}

\Needspace{5\baselineskip}
\begin{example}\label{ex:A-reading}
Comparing coefficients in $\log A=\frac{z(1+A)}{2}$ gives
\[
 A(z)=1+z+z^2+\frac{7}{6}z^3+O(z^4).
\]
Thus $f_1=1$, $f_2=1$, and $f_3=2$. By the normalization in~\eqref{eq:Adef},
$f_4=3!\,[z^3]A=7$, so there are exactly seven feasible trees on the
vertex set $\{1,2,3,4\}$.
\end{example}

\section{Root capacity by rank}
\label{sec:G-refinements}
Having determined $A$, we return to $G$ to count feasible trees by
root rank and capacity. By~\eqref{eq:FGparts} and
Remark~\ref{lem:feasibility}, the contribution to $[z^n]G$ of a tree
$T\in\Fn$ rooted at a vertex $r$ of rank $j$ is the factor for rank $j$
in~\eqref{eq:Fdef} multiplied by $t^{\kappa_T(r)}$.
To determine the three series $a,b,d$ in~\eqref{eq:rationalG}, set
\begin{equation}\label{eq:Edef}
 E(x;z)=\exp\bigl(z(1+A(z))x\bigr),\qquad 0\le x\le1.
\end{equation}

\begin{theorem}\label{thm:refinement}
With $A=A(z)$ and $E=E(x;z)$, the three series in~\eqref{eq:rationalG} are
\begin{align}
 a(x)&=\frac{z(A^2-E)}{A+E},\label{eq:a-explicit}\\
 b(x)&=\frac{zE(A+1)^2}{(A+E)^2},\label{eq:b-explicit}\\
 d(x)&=\frac{E-1}{A+E}.\label{eq:d-explicit}
\end{align}
\end{theorem}

\begin{proof}
By~\eqref{eq:rho} and the definition of $E$ in~\eqref{eq:Edef},
\[
 \frac{z+Cd}{1-d}=zE.
\]
Since $C=zA$, cancelling $z$ gives $1+Ad=E(1-d)$, which yields~\eqref{eq:d-explicit}.
Substitution in the two identities of Corollary~\ref{cor:invariants} gives
\[
 a=C-z-Cd=\frac{z(A^2-E)}{A+E},\qquad
 b=(1-d)(z+Cd)=\frac{zE(A+1)^2}{(A+E)^2}.
\]
The denominator $A+E$ has constant term $2$ in $z$, so it is invertible as a formal power series.
\end{proof}

\begin{example}\label{ex:G-reading}
The two trees in $\mathcal F_3$ are $T_1$ and $T_2$ in Figure~\ref{fig:capacity-n3}.
Substituting the three series from Theorem~\ref{thm:refinement} into~\eqref{eq:rationalG} and taking the coefficient of $z^3$ gives
\begin{align*}
 [z^3]G(x,t;z)
 &=1-x+x(1-x)t+x^2t^3\\
 &=\frac{(1-x)^2}{2!}\,2
   +x(1-x)(1+t)
   +\frac{x^2}{2!}\,2t^3.
\end{align*}
In the decomposition by rank, the three summands correspond to roots
$1$, $2$, and $3$, in that order.
For root $2$, the rank factor is $x(1-x)$, and the factor $1+t$ comes
from the two trees in Figure~\ref{fig:middle-root}.

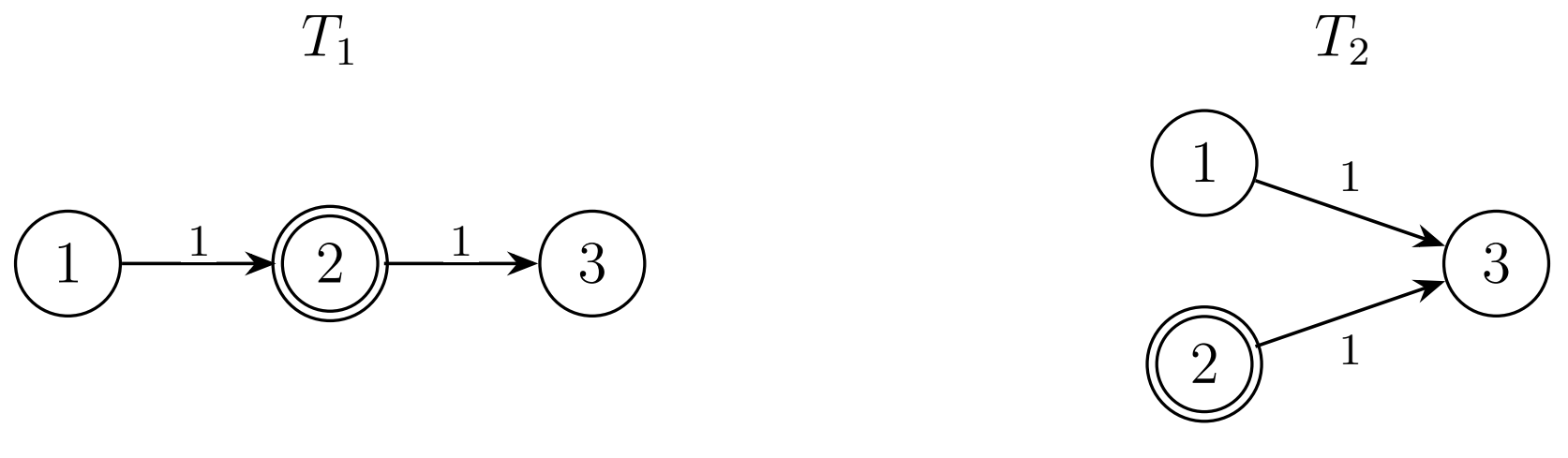
\begin{figure}[H]
\centering
\begin{tikzpicture}[
  vertex/.style={circle,draw,line width=.62pt,minimum size=7.4mm,inner sep=0pt,font=\large},
  root/.style={vertex,double,double distance=1.3pt},
  edge/.style={-{Stealth[length=2.2mm,width=1.7mm]},line width=.68pt},
  weight/.style={fill=white,inner sep=1.2pt,font=\small}]
\begin{scope}
  \node at (0,1.65) {\large $T_1$};
  \node[vertex] (d) at (-1.85,.07) {$1$};
  \node[root]   (e) at (0,.07) {$2$};
  \node[vertex] (f) at (1.85,.07) {$3$};
  \draw[edge] (d)--node[weight,above] {$1$}(e);
  \draw[edge] (e)--node[weight,above] {$1$}(f);
\end{scope}
\begin{scope}[xshift=7.15cm]
  \node at (0,1.65) {\large $T_2$};
  \node[vertex] (a) at (-.98,.78) {$1$};
  \node[root]   (b) at (-.98,-.64) {$2$};
  \node[vertex] (c) at (1.08,.07) {$3$};
  \draw[edge] (a)--node[weight,midway,above=3pt] {$1$}(c);
  \draw[edge] (b)--node[weight,midway,below=3pt] {$1$}(c);
\end{scope}
\end{tikzpicture}
\caption{The displayed weights attain the capacity at vertex $2$ in each tree.}
\label{fig:middle-root}
\end{figure}

On the path $T_1$, feasibility at vertex $1$ forces $w(12)=1$.
Thus the divergence at vertex $2$ is $w(23)-1$, whose minimum is $0$,
attained when $w(23)=1$. Hence $\kappa_{T_1}(2)=1$. The path has two
feasible weightings and contributes a single $t$.

In $T_2$, the only edge at vertex $2$ is the outgoing edge $2\to3$.
Feasibility forces its weight to be $1$, so the divergence at $2$ is $1$.
Thus $\kappa_{T_2}(2)=0$, and $T_2$ contributes the constant term $1$.
\end{example}

\begin{corollary}\label{cor:spec}
Recall that $C=zA$. The endpoint value of $d$ is
$d(1)=1-\frac{1}{A}$. For feasible trees with a distinguished root,
integrating $G(x,t;z)$ from $x=0$ to $x=1$ gives
\begin{equation}\label{eq:allroots}
 \sum_{n\ge1}\frac{z^n}{n!}\sum_{T\in\Fn}\sum_{r\in V(T)} t^{\kappa_T(r)}
 =\frac{C-z}{2}-\log(1-t\,d(1)).
\end{equation}
For feasible trees rooted at their largest vertex, setting $x=1$ in
$G(x,t;z)$ gives
\begin{equation}\label{eq:maxroot}
 \sum_{n\ge1}\frac{z^n}{(n-1)!}\sum_{T\in\Fn} t^{\kappa_T(n)}
 =\frac{zt}{1-t\,d(1)}.
\end{equation}
In particular, for $k\ge1$,
\begin{equation}\label{eq:capacity-k-all}
 \sum_{n\ge1}\frac{z^n}{n!}
 \bigl|\{(T,r):T\in\Fn,\ r\in V(T),\ \kappa_T(r)=k\}\bigr|
 =\frac{d(1)^k}{k},
\end{equation}
and
\begin{equation}\label{eq:capacity-k-max}
 \sum_{n\ge1}\frac{z^n}{(n-1)!}
 \bigl|\{T\in\Fn:\ \kappa_T(n)=k\}\bigr|
 =z\,d(1)^{k-1}.
\end{equation}
\end{corollary}

\begin{proof}
For each root rank $j$, the beta integral gives
\[
 \int_0^1\frac{x^{j-1}(1-x)^{n-j}}{(j-1)!(n-j)!}\,dx=\frac{1}{n!}.
\]
Thus integrating the rank factor in~\eqref{eq:Fdef} gives each rooted pair $(T,r)$ the factor $\frac{1}{n!}$, so the left side of~\eqref{eq:allroots} is $\int_0^1 G(x,t;z)\,dx$.
By~\eqref{eq:Edef}, $E(0;z)=1$. Equation~\eqref{eq:Aeq} gives
\[
 E(1;z)=\exp\bigl(z(1+A)\bigr)=A^2.
\]
Hence~\eqref{eq:d-explicit} yields
\[
 d(0)=0,\qquad d(1)=\frac{A^2-1}{A+A^2}=1-\frac{1}{A}.
\]
Using $d'=b$ from~\eqref{eq:abd}, we obtain
\[
 \frac{d}{dx}\bigl[-\log(1-d(x)t)\bigr]
 =\frac{d'(x)t}{1-d(x)t}=\frac{b(x)t}{1-d(x)t}.
\]
Integrating coefficientwise in $z$ gives
\[
 \int_0^1\frac{b(x)t}{1-d(x)t}\,dx
 =\left[-\log(1-d(x)t)\right]_{x=0}^{x=1}
 =-\log(1-t\,d(1)).
\]
At $t=1$, Lemma~\ref{lem:C} gives $G(x,1;z)=C$.
Since $\log A=\frac{C+z}{2}$ by~\eqref{eq:Aeq}, we obtain
\[
 \int_0^1 a(x)\,dx=C+\log(1-d(1))=C-\log A=\frac{C-z}{2}.
\]
This proves~\eqref{eq:allroots}.

With $E(1;z)=A^2$, equations~\eqref{eq:a-explicit}
and~\eqref{eq:b-explicit} give $a(1)=0$ and $b(1)=z$.
Evaluating~\eqref{eq:rationalG} at $x=1$ with these endpoint values gives
$G(1,t;z)=\frac{zt}{1-t\,d(1)}$. At $x=1$, only rank $j=n$ contributes in~\eqref{eq:Fdef}, with factor $\frac{1}{(n-1)!}$. This gives~\eqref{eq:maxroot}. Extracting coefficients of $t$ proves the remaining identities.
\end{proof}

\begin{example}\label{ex:capacity-specializations}
From the expansion of $A$ in Example~\ref{ex:A-reading},
\[
 C=z+z^2+z^3+O(z^4),\qquad d(1)=z+\frac{z^3}{6}+O(z^4).
\]
The two specializations in Corollary~\ref{cor:spec} give
\[
 3![z^3]\left(\frac{C-z}{2}-\log(1-t\,d(1))\right)=3+t+2t^3,\qquad
 2![z^3]\frac{zt}{1-t\,d(1)}=2t^3.
\]
Of the six rooted pairs $(T,r)$ counted by $3+t+2t^3$, three have root capacity $0$, one has capacity $1$, and two have capacity $3$. When each tree is rooted at its largest vertex, $3$, both have root capacity $3$, so each contributes $t^3$ to the second polynomial.

The weights in Figure~\ref{fig:middle-root} attain the capacity at vertex $2$. At vertex $3$ of the path $T_1$, the slack is $1+w(23)$. Since $w(12)=1$, the bound at vertex $2$ gives $w(23)-1\le1$. Taking $w(23)=2$ therefore gives $\kappa_{T_1}(3)=3$.

\begin{table}[H]
\centering
\begin{tabular}{c|ccc}
 & $\kappa_T(1)$ & $\kappa_T(2)$ & $\kappa_T(3)$\\
\hline
$T_1$ & $0$ & $1$ & $3$\\
$T_2$ & $0$ & $0$ & $3$
\end{tabular}
\caption{Root capacities of the two trees in $\mathcal F_3$.}
\label{tab:root-capacities}
\end{table}

Summing $t^{\kappa_T(r)}$ over all six entries of Table~\ref{tab:root-capacities} recovers the first polynomial, and summing over the last column recovers the second.
\end{example}

Call a vertex of a feasible tree \emph{free} if its capacity is positive and
\emph{forced} if its capacity is zero. At a forced vertex every feasible weighting
has divergence $1$.

\begin{corollary}\label{cor:free}
For $n\ge2$ and each rank $1\le j\le n$, the numbers of forced roots
of rank $j$ and free roots of rank $n+1-j$ satisfy
\begin{equation}\label{eq:reflection}
 \bigl|\{T\in\Fn:\kappa_T(j)=0\}\bigr|
 =\bigl|\{T\in\Fn:\kappa_T(n+1-j)>0\}\bigr|.
\end{equation}
Consequently, exactly half of the pairs $(T,r)$ with $T\in\Fn$ and $r\in V(T)$
have positive root capacity. Thus
\begin{equation}\label{eq:half}
 \sum_{T\in\Fn}|\{r\in V(T):\kappa_T(r)>0\}|=\frac{nf_n}{2}.
\end{equation}
\end{corollary}

\begin{proof}
By~\eqref{eq:Edef} and the endpoint value obtained in the proof of
Corollary~\ref{cor:spec},
\[
 E(x;z)E(1-x;z)=\exp\bigl(z(1+A)\bigr)=E(1;z)=A^2.
\]
Writing $E=E(x;z)$, we have $E(1-x;z)=\frac{A^2}{E}$. Substituting into
\eqref{eq:a-explicit} gives the two terms
\[
 a(x)=\frac{z(A^2-E)}{A+E},\qquad
 a(1-x)=\frac{zA^2(E-1)}{AE+A^2}
       =\frac{zA(E-1)}{A+E}.
\]
Adding them gives
\[
 a(x)+a(1-x)
 =\frac{z\bigl(A^2-E+A(E-1)\bigr)}{A+E}
 =z(A-1)=C-z.
\]
By~\eqref{eq:Fdef} and~\eqref{eq:FGparts}, $a(x)=G(x,0;z)$ counts feasible trees whose root has capacity zero. For $n\ge2$, expand the coefficient of $z^n$ on both sides in the rank basis from~\eqref{eq:Fdef},
\[
 \frac{x^{j-1}(1-x)^{n-j}}{(j-1)!(n-j)!},\qquad 1\le j\le n.
\]
Replacing $x$ by $1-x$ interchanges the basis terms for ranks $j$ and $n+1-j$. By~\eqref{eq:Cdef}, $[z^n](C-z)=\frac{f_n}{(n-1)!}$. Since the basis terms sum to $\frac{1}{(n-1)!}$, each has coefficient $f_n$ on the right. Comparing coefficients gives
\[
 \bigl|\{T\in\Fn:\kappa_T(j)=0\}\bigr|
 +\bigl|\{T\in\Fn:\kappa_T(n+1-j)=0\}\bigr|=f_n.
\]
Since capacities are nonnegative, this is equivalent to~\eqref{eq:reflection}. Summing over the ranks $j$ and counting the corresponding vertices $r$ gives
\[
 2\sum_{T\in\Fn}\bigl|\{r\in V(T):\kappa_T(r)=0\}\bigr|=nf_n.
\]
There are $nf_n$ pairs $(T,r)$ in total, so~\eqref{eq:half} follows.
\end{proof}

\begin{corollary}\label{cor:total}
For every $n\ge1$,
\begin{equation}\label{eq:totalcap}
 \sum_{T\in\Fn}\sum_{r\in V(T)}\kappa_T(r)=f_{n+1}.
\end{equation}
\end{corollary}

\begin{proof}
Differentiate~\eqref{eq:allroots} with respect to $t$ and set $t=1$.
Each rooted pair $(T,r)$ contributes $\kappa_T(r)$, since the derivative
of $t^{\kappa_T(r)}$ at $t=1$ is $\kappa_T(r)$. Hence
\[
 \sum_{n\ge1}\frac{z^n}{n!}\sum_{T\in\Fn}\sum_{r\in V(T)}\kappa_T(r)
 =\left.\frac{d(1)}{1-t\,d(1)}\right|_{t=1}
 =\frac{d(1)}{1-d(1)}.
\]
Since $1-d(1)=\frac{1}{A}$, the last quotient equals $A-1$. By~\eqref{eq:Adef}
and the fact that $f_1=1$,
\[
 A(z)-1
 =\sum_{m\ge2}f_m\frac{z^{m-1}}{(m-1)!}
 =\sum_{n\ge1}f_{n+1}\frac{z^n}{n!}.
\]
Comparing coefficients of $\frac{z^n}{n!}$ gives~\eqref{eq:totalcap}.
\end{proof}

\begin{example}\label{ex:free-total}
At ranks $1,2,3$, Table~\ref{tab:root-capacities} gives $2,1,0$ forced roots and $0,1,2$ free roots, respectively. Reversing the ranks interchanges the two lists, as in Corollary~\ref{cor:free}. Three of the six rooted pairs are therefore free. The two rows have total capacities $4$ and $3$, giving $4+3=7=f_4$ as in Corollary~\ref{cor:total}.
\end{example}

Theorem~\ref{thm:main} suggests looking for a direct bijection between $\Fn$ and the alternating trees on the vertex set $\{1,\ldots,n\}$.

\section*{Acknowledgement}
{\raggedright
This work was supported by the Natural Science Foundation of
Tianjin Municipality (Grant No.\ \mbox{25JCQNJC00250}) and the National Natural
Science Foundation of China (Grant No.~12001398).\par}

\section*{Declaration of AI Use}
During the preparation of this manuscript, ChatGPT (OpenAI) was used
as an auxiliary tool for exploratory discussions of possible proof strategies
and for English language editing. The authors reviewed the suggestions
incorporated into the manuscript and take full responsibility for all mathematical
results, proofs, and the final text.

\Needspace{5\baselineskip}
\begin{flushleft}
\small
Department of Mathematics, Tianjin University of Finance and Economics,
Tianjin 300222, P.R. China\par\smallskip
\begin{tabular}{@{}l@{\quad}l@{}}
Xuanying Han & \href{mailto:hanxuanying@stu.tjufe.edu.cn}{\texttt{hanxuanying@stu.tjufe.edu.cn}} \\
Sophie C.C. Sun & \href{mailto:sophiesun@tjufe.edu.cn}{\texttt{sophiesun@tjufe.edu.cn}}
\end{tabular}
\end{flushleft}

\end{document}